\documentclass[sn-mathphys-num]{sn-jnl}

\usepackage{graphicx}%
\usepackage{multirow}%
\usepackage{amsmath,amssymb,amsfonts}%
\usepackage{amsthm}%
\usepackage{mathrsfs}%
\usepackage[title]{appendix}%
\usepackage{xcolor}%
\usepackage{textcomp}%
\usepackage{manyfoot}%
\usepackage{booktabs}%
\usepackage{algorithm}%
\usepackage{algorithmicx}%
\usepackage{algpseudocode}%
\usepackage{listings}%
\usepackage{comment}
\usepackage[colorinlistoftodos]{todonotes}
\usepackage{tabularx}
\usepackage{ragged2e}
\usepackage{makecell}
\usepackage{subfigure}
\usepackage{threeparttable}

\theoremstyle{thmstyleone}%
\theoremstyle{thmstyletwo}%

\theoremstyle{thmstylethree}%

\begin{document}

\title[Article Title]{Deterministic Global Optimization over trained Kolmogorov Arnold Networks}


\author[1]{\fnm{Tanuj} \sur{Karia}}

\author[1]{\fnm{Giacomo} \sur{Lastrucci}}

\author*[1]{\fnm{Artur M.} \sur{Schweidtmann}}\email{a.schweidtmann@tudelft.nl}

\affil[1]{\orgdiv{Department of Chemical Engineering}, \orgname{Delft University of Tecnology}, \orgaddress{\street{Van der Maasweg 9}, \city{Delft}, \postcode{2629 HZ}, 
\country{The Netherlands}}}




\abstract{To address the challenge of tractability for optimizing mathematical models in science and engineering, surrogate models are often employed. Recently, a new class of machine learning models named Kolmogorov Arnold Networks (KANs) have been proposed. It was reported that KANs can approximate a given input/output relationship with a high level of accuracy, requiring significantly fewer parameters than multilayer perceptrons. Hence, we aim to assess the suitability of deterministic global optimization of trained KANs by proposing their Mixed-Integer Nonlinear Programming (MINLP) formulation. We conduct extensive computational experiments for different KAN architectures. Additionally, we propose alternative convex hull reformulation, local support and redundant constraints for the formulation aimed at improving the effectiveness of the MINLP formulation of the KAN. KANs demonstrate high accuracy while requiring relatively modest computational effort to optimize them, particularly for cases with less than five inputs or outputs. For cases with higher inputs or outputs, carefully considering the KAN architecture during training may improve its effectiveness while optimizing over a trained KAN. Overall, we observe that KANs offer a promising alternative as surrogate models for deterministic global optimization.}

\keywords{Deterministic global optimization, Mixed-Integer Nonlinear Programming, Kolmogorov Arnold Networks, Machine Learning}



\maketitle
\section{Introduction}\label{sec1:intro}
Surrogate models are frequently employed for mathematical optimization in various disciplines of science and engineering \citep{Fischetti2019,Abbasi2020,Schweidtmann2020,Misener2023}. 
Often, when solving optimization problems, one may encounter either an objective or a set of constraints, which makes the problem challenging to solve.
In such cases, surrogate models are often trained by querying samples from the objective or set of constraints to approximate the underlying input-output relationship exhibited by the objective or constraint.
Typically, the surrogate model has a simpler functional form and requires less computational effort to optimize, possibly at some loss of accuracy.
\par
Various kinds of surrogate models have been used for mathematical optimization, including linear regression models \citep{Dowling2014}, polynomial regression models \citep{Palmer2002}, and other algebraic models \citep{Ma2022,Ma2022a} derived via sparse regression techniques \citep{Cozad2018}. 
Other surrogate model forms \citep{Schweidtmann2020} include  convex region linear surrogate models \citep{Zhang2016}, piecewise polynomial \citep{Grimstad2020} and spline functions \citep{Grimstad2016, Grimstad2018}, Gaussian processes \citep{Schweidtmann2021}, support vector machines \citep{Beykal2020,Schweidtmann2022}, ensemble tree models \citep{Misic2020,Mistry2021,Thebelt2021,Ammari2023}, graph neural networks \citep{McDonald2024,Zhang2024}, and artificial neural networks or multilayer perceptrons (MLPs) \citep{Schweidtmann2019,Grimstad2019,Anderson2020,Tsay2021,Wilhelm2023,Wang2023}.
\par
Besides, the development of various modeling frameworks has greatly facilitated the integration of surrogate models into optimization formulations. 
For embedding \textsf{ReLU} networks into optimization models, several tools namely \texttt{JANOS} \citep{Bergman2022}, \texttt{optiCL} \citep{Maragno2023}, \texttt{reluMIP} \citep{Lueg2021} have been developed.
Catering to more general machine learning models, \texttt{OMLT} \citep{Ceccon2022} is capable of incorporating MLPs, tree ensemble models and graph neural networks \citep{McDonald2024,Zhang2024} into optimization models formulated in \texttt{Pyomo} \citep{Bynum2021}.
\texttt{OMLT} allows for seamless switching between various formulations for the same machine learning model.
\texttt{MeLON} \citep{Schweidtmann2019,Schweidtmann2021} is capable of handling Gaussian processes, MLPs and, support vector machines. \texttt{MeLON} interfaces with C++ for the global solver \texttt{MAiNGO} \citep{Bongartz2018}.
Recently, more solver-specific interfaces also have been developed for \texttt{GUROBI} \citep{Gurobi2023} and \texttt{SCIP} \citep{Turner2023}.
Additionally, a package named MathOptAI has been published for embedding machine learning models in JuMP \citep{MathOptAI}.
For a more detailed review of modeling frameworks for integrating surrogate models into optimization formulations, the reader is referred to the article by \citet{LopezFlores2024}.
\par
The popularity of deep learning, in combination with the availability of a wide plethora of tools for integrating MLPs, has made MLPs a popular choice as surrogate models. 
The most popular activation functions employed for the training of MLPs are the \textsf{tanh} and \textsf{ReLU} activations.
Significant work has been done in the literature to facilitate the optimization of deeper (more layers) and wider (more neurons) MLPs with these activation functions.
\citet{Schweidtmann2019} proposed convex and concave envelopes for the \textsf{tanh} function and employed a factorable reduced space formulation of the trained MLP to exploit lower dimensionality to effectively optimize ANNs with \textsf{tanh} activation.
Alternatively, using \textsf{ReLU} activation function results in a mixed-integer linear formulation \citep{Grimstad2019} because of the use of big-M method to model the \textsf{max} function.
Stronger \citep{Anderson2020} and more effective partition-based \citep{Tsay2021} mixed-integer linear formulations have further improved the tractability of optimizing over trained \textsf{ReLU} networks.
Recently, a complementarity constraint-based formulation for modelling the \textsf{ReLU} activation leading to a nonlinear formulation was also proposed \citep{Yang2022}.
Finally, convex and concave relaxations for other activation functions such as sigmoid, GeLU, and SiLU have been developed \citep{Wilhelm2023,carrasco2024tighteningconvexrelaxationstrained} to facilitate efficient optimization of MLPs with these activation functions. 
Generally, it has been reported in the literature that the larger the ANN, both in terms of number of neurons and layers, the more challenging it becomes to solve the resulting optimization model \citep{Schweidtmann2019,Gurobi2023}.
\par
Recently, a new class of machine-learning models, Kolmogorov-Arnold Networks (KANs), have been proposed \citep{Liu2024}.
\citet{Liu2024} demonstrated the superiority of KANs over MLPs for different machine-learning tasks in terms of the number of parameters required to achieve a similar level of accuracy at the expense of slower training times for KANs.
Herein, we investigate the suitability of KANs as surrogate models in mathematical optimization, particularly considering that KANs are expected to require smaller architectures than MLPs to approximate a given input/output relationship.
KANs are network-type models similar to multi-layer perceptrons with two features distinguishing them from MLPs:
(i) the activation lies on the edge (weight) connecting two neurons in a KAN rather than on the neuron itself for a MLP, and (ii) each activation is learnable via B-spline functions.

Thus, we consider optimization problems of the form:
\begin{align}
    \min_{\mathbf{x}} & \quad f(\mathbf{x}, \textsf{KAN}(\mathbf{x})) \label{eq:objective}\\
    \text{s.t.} & \quad \mathbf{h}(\mathbf{x},\textsf{KAN}(\mathbf{x})) = 0 \label{eq:equality} \\
    & \quad \mathbf{g}(\mathbf{x},\textsf{KAN}(\mathbf{x})) \leq 0 \label{eq:inequality} \\
    & \quad \mathbf{x} \in \mathcal{X} \subset \mathbb{R}^{n_0} \label{eq:bounds}
\end{align}
Here, $f$ is the objective function, $\mathbf{h}$ and, $\mathbf{g}$ represents the set of equality and inequality constraints, respectively, describing the formulation of a trained KAN; $\textsf{KAN}(\mathbf{x}) \in Y \subset \mathbb{R}^{n_L}$ is the output of the KAN at $\mathbf{x}$. 
\par
We propose a MINLP formulation for trained Kolmogorov-Arnold Networks, which is based on the Mixed-Integer Quadratically Constrained Programming (MIQCP) formulation proposed for B-splines by \citet{Grimstad2018}. 
Further, we propose enhancements to exploit some properties of the KANs to improve the effectiveness of the proposed MINLP formulation. 
The proposed formulation is implemented using \texttt{Pyomo} \citep{Bynum2021} and tested on various test functions for optimization of varying dimensionality. 
Different architectures of KANs are trained and optimized using global MINLP solver \texttt{SCIP} \citep{Bolusani2024} to conduct a rigorous examination of whether KANs are suitable as surrogate models. 
MLPs of varying architectures are also trained, and the computational effort needed to solve both MLPs and KANs to global optimality with \texttt{SCIP} is compared.
\par
The rest of the paper is as follows: Section 2 provides a brief introduction describing KANs for the reader. 
Section 3 outlines the methodology in which the MINLP formulation for KANs and the enhancements to the formulation are explained.
The implementation details in \texttt{Pyomo} and details of the computational experiments are also provided in Section 3. 
Section 4 provides the results and discussion on optimizing the trained KANs. 
Finally, we conclude and provide recommendations for future research in Section 5 of the paper.

\section{Background}\label{sec2:background}
In this section, we provide the requisite background on Kolmogorov Arnold Networks (KANs) and B-spline functions which act as activations in a KAN.
\subsection{Kolmogorov Arnold Networks}
Kolmogorov Arnold Networks are based on Kolmogorov-Arnold representation theorem \citep{Kolmogorov1957,Givental2009} which establishes that for a given vector $\mathbf{x}$ of $n_0$ variables, a continuous function $f(\mathbf{x})$ on a bounded domain $\mathcal{X} \in \mathbb{R}^{n_0}$ such that $f:\mathcal{X} \rightarrow \mathbb{R}$, can be expressed with ($2n_0 + 1$) additions of a function $\Phi$, which is a finite composition of $n_0$ univariate functions $\phi$.
Specifically,
\begin{equation}
    \label{eq:kol-arn-rep-theorem}
    f(\mathbf{x}) = f(x_1, \ldots, x_{n_0}) = \sum_{j=1}^{2n_{0}+1} \Phi_j \left( \sum_{i=1}^{n_0} \phi_{j,i}(x_i) \right)
\end{equation}
In Equation \eqref{eq:kol-arn-rep-theorem}, $\phi_{j,i}(x_i):[x_i^{\rm L}, x_i^{\rm U}] \rightarrow \mathbb{R}$ and $\Phi_j:\mathbb{R} \rightarrow \mathbb{R}$ represents the composition of univariate functions $\phi_{j,i}$.
Kolmogorov-Arnold representation theorem thereby opens up the possibility of learning any function in $n_0$-dimensions using a polynomial number of 1D functions.
Hence, a network based on the Kolmogorov-Arnold representation theorem comprises one input layer of size $n_0$, one output layer with output $f(\mathbf{x})$ and, one hidden layer with ($2n_0+1$) neurons.
However, the 1D functions can be non-smooth or fractal, rendering them unlearnable in practice \citep{Girosi1989,Poggio2020}.
\begin{figure}[H]
    \centering
    \includegraphics[width=\textwidth]{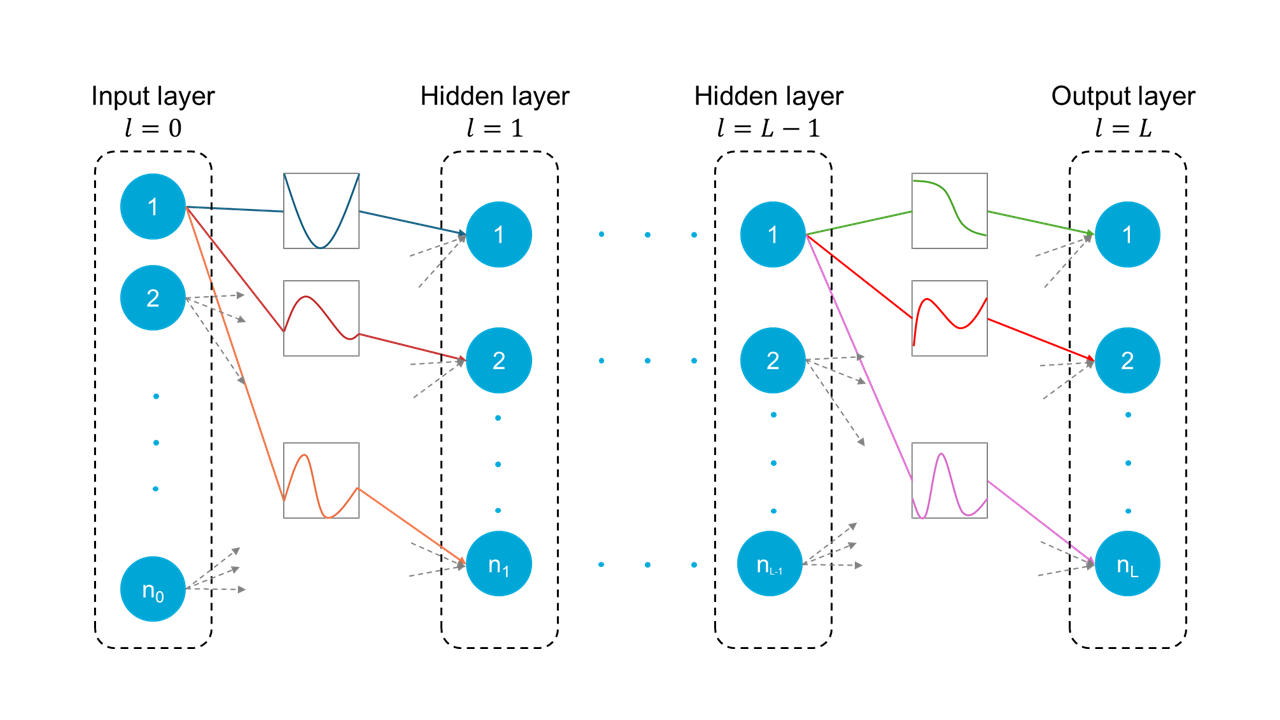}
    \caption{Representation of a Kolmogorov-Arnold Network as a directed acyclic graph. Dashed gray line on a neuron represents the edge with the activation connecting two neurons. The network is fully connected.}
    \label{fig:kan-dag}
\end{figure}
To address this issue, \citet{Liu2024} proposed generalizing Equation~\eqref{eq:kol-arn-rep-theorem} by considering an arbitrarily large number of hidden layers (depth of the network) and neurons (width) of the hidden layers 
and applied modern backpropagation techniques to train them. 
A generalised KAN is illustrated as a directed acyclic graph in Figure~\ref{fig:kan-dag}.

\par
Assuming a set of $N_{\rm data}$ input-output pairs $(\mathbf{x}, \mathbf{y})$ is available, the aim is to approximate the relationship $\mathbf{y} \approx f(\mathbf{x})$ using a KAN.
Considering a KAN of arbitrary width and depth comprising $L+1$ layers, including the input ($l=0$) and output ($l=L$) layers, such that the number of neurons in each layer $l$ is defined by the vector $\mathbf{N} := [n_0, n_1, \ldots, n_l, \ldots, n_{L-1}, n_{\rm L}]$ of length $L+1$.
As seen in Figure~\ref{fig:kan-dag}, the univariate activation function for a KAN ($\phi$) lies on the edge connecting two neurons.
The univariate activation function ($\phi$) connects a neuron $i$ in layer $l$ to a neuron $j$ in the next layer $l+1$.
Hence, the univariate activation function ($\phi$) is indexed as
\begin{equation}
    \label{eq: activation-index}
    \phi_{l,j,i}(x_{l,i}),\, l=0,\ldots,\,L-1, \, j=1,\ldots, \, n_{l+1},\, \text{and } i=1,\ldots,\, n_{l}
\end{equation}
In Equation~\eqref{eq: activation-index}, $x_{l,i}$ refers to the input $x$ to the activation function in layer $l$ from neuron $i$.
Based on the Kolmogorov representation theorem introduced earlier, the input ($x_{l+1,j}$) to the neuron $j$ in the next layer $l+1$ is a sum of all the incoming outputs of the activations from layer $l$ to the neuron $j$,
\begin{equation}
    \label{eq: output-scalar-form}
    x_{l+1,j} = \sum_{i=1}^{n_l} \phi_{l,j,i}(x_{l,i}), \, j=1,\ldots,\,n_{l+1}
\end{equation}
Equivalently, Equation~\eqref{eq: output-scalar-form} can be represented as matrices,
\begin{equation}
    \begin{pmatrix}
        x_{l+1,1} \\
        \vdots \\
        x_{l+1,n_{l+1}}
    \end{pmatrix}
    =
    \begin{pmatrix}
        \phi_{l,1,1}(x_{l,1}) + \phi_{l,1,2}(x_{l,2}) + \cdots + \phi_{l,1,2}(x_{l,n_l}) \\
        \phi_{l,2,1}(x_{l,1}) + \phi_{l,2,2}(x_{l,2}) + \cdots + \phi_{l,2,2}(x_{l,n_l}) \\
        \vdots                                                                              \\
        \phi_{l,n_{l+1},1}(x_{l,1}) + \phi_{l,n_{l+1},2}(x_{l,2}) + \cdots + \phi_{l,n_{l+1},2}(x_{l,n_l})
    \end{pmatrix}
    \label{eq: output-matrix-form}
\end{equation}
Equation~\eqref{eq: output-matrix-form} can be equivalently represented as,
\begin{equation}
    \begin{pmatrix}
        x_{l+1,1} \\
        \vdots \\
        x_{l+1,n_{l+1}}
    \end{pmatrix}
    =
    \begin{pmatrix}
        \phi_{l,1,1}(\cdot) & \phi_{l,1,2}(\cdot) & \cdots & \phi_{l,1,2}(\cdot) \\
        \phi_{l,2,1}(\cdot) & \phi_{l,2,2}(\cdot) & \cdots & \phi_{l,2,2}(\cdot) \\
        \vdots              & \vdots              & \ddots & \vdots               \\
        \phi_{l,n_{l+1},1}(\cdot) & \phi_{l,n_{l+1},2}(\cdot) & \cdots & \phi_{l,n_{l+1},2}(\cdot)
    \end{pmatrix}
    \begin{pmatrix}
        x_{l,1} \\
        \vdots \\
        x_{l,n_{l}}
    \end{pmatrix}
    \label{eq: output-matrix-form-form-function-rep}
\end{equation}
It is important to note that in Equation~\eqref{eq: output-matrix-form-function-rep}, the second matrix on the right-hand side of the equation is not multiplied to the first matrix on the right-hand side of the equation.
Rather, this shows all elements in a given row in the first matrix take in the element from the corresponding column in the second matrix as input to the activation $\phi$.
Compactly, Equation~\eqref{eq: output-matrix-form-form-function-rep} is represented as,
\begin{equation}
    \mathbf{x}_{l+1} = \Phi_l(\mathbf{x}_{l})
    \label{eq: output-matrix-form-compact}
\end{equation}
In Equation~\eqref{eq: output-matrix-form-compact}, $\Phi_l$ is a ($n_{l+1} \times n_l$) matrix representing the composition of univariate functions $\phi$, which is a sum of all activations from neuron $i$ in layer $l$ to neuron $j$ in layer $l+1$. $n_l$ represents the number of inputs to the layer $l$ and $n_{l+1}$ represents the number of outputs from layer $l$ or the number of inputs to the layer $l+1$.
Hence, the predicted output from the KAN ($\mathbf{y}^{\rm pred}$) can be expressed compactly as
\begin{equation}
    \label{eq: kan-compact-form}
    \mathbf{y}^{\rm pred} \equiv \mathbf{x}_{L} = \Phi_{L-1}(\cdots(\Phi_1(\Phi_0(\mathbf{x}_0))))
\end{equation}
, where $\mathbf{x}$ is equivalent to $\mathbf{x}_0$ representing the inputs to the network.
\par
As described in Figure~\ref{fig:kan-dag}, in a KAN, each activation function ($\phi$) is trainable; thereby, each activation in the network is distinct.
Each univariate activation function $\phi$ is defined as:
\begin{align}
    \phi(x_{l,i}) & = w^{\rm b} b(x_{l,i}) + w^{\rm s}s(x_{l,i}) \label{eq:kan-function} \\
    b(x_{l,i}) & = \text{silu}(x_{l,i}) = \frac{x_{l,i}}{1 + \exp{(-x_{l,i})}} \label{eq:kan-function-base} \\
    s(x_{l,i}) & = \sum_{g=0}^{G+k-1} c_g B_{g,k}(x_{l,i}) \label{eq:kan-function-spline}
\end{align}
In Equations~\eqref{eq:kan-function} - \eqref{eq:kan-function-spline}, the indices $l,j,i$ are dropped for simplicity but these equations are applicable for all $l=0,\ldots,\,L-1, \, j=1,\ldots, \, n_{l+1},\, \text{and } i=1,\ldots,\, n_{l}$. 
A KAN activation is a weighted sum of the univariate spline activation function ($s(x_{l,i})$) and a univariate base activation function ($b(x_{l,i})$) as described in Equation~\eqref{eq:kan-function}, where $w^{\rm s}$ and $w^{\rm b}$ represent its weights, respectively.
In principle, any activation function can be chosen as the base activation function, but \citet{Liu2024} propose the use of \texttt{SiLU} activation in Equation~\eqref{eq:kan-function-base} to ensure stability during training.
Each activation ($\phi$) in a KAN is distinct due to the presence of the univariate spline activation function ($s(x_{l,i})$).
The univariate spline function is expressed as a sum-product of $k$\textsuperscript{th}-order B-spline basis functions ($B_{g,k}$) and their respective coefficients ($c_g$) over the set $g \in \mathcal{G}$.
The set $\mathcal{G}$ is a non-decreasing sequence of real numbers called knots.
The number of elements in the set $\mathcal{G}$ depends on the number of grid points $G$ and the order $k$ of the B-spline chosen for training by the user, thereby acting as hyperparameters.
The bias-variance tradeoff for KANs is determined by the number of intervals $G$ used to train KANs \citep{Liu2024}.
\par
For a given input-output pair $(\mathbf{x}, \mathbf{y})$ and number of grid points $G$, the input $\mathbf{x}_{0}$ corresponds to $\mathbf{x}$. 
The objective is to ensure the network output $\mathbf{x}_{L} \equiv \mathbf{y}^{\rm pred}$ closely matches the value of $\mathbf{y}$.
This is done by adjusting the values of $c_g$, $w^{\rm b}$, and $w^{\rm s}$, thus representing the trainable parameters of a KAN.
During the training of a KAN, $B_{g,k}$ are fixed parameters, and they take values based on the value of the input ($x_{0,i}$) to the network.
The value of $B_{g,k}$ is determined based on the de Boor algorithm \citep{Boor1977} and is described in Equations~\eqref{eq:0-basis-function} and \eqref{eq:other-basis-functions}.
\begin{align}
B_{g,0}(x_{l,i};\mathbf{t}) &= 
\begin{cases}
1, & t_g \leq x_{l,i} < t_{g+1}, \\
0, & \text{otherwise}.
\end{cases} \label{eq:0-basis-function} \\
B_{g,d}(x_{l,i};\mathbf{t}) &= \frac{x_{l,i} - t_g}{t_{g+d} - t_g} B_{g,d-1}(x_{l,i};\mathbf{t}) + \frac{t_{g+d+1} - x_{l,i}}{t_{g+d+1} - t_{g+1}} B_{g+1,d-1}(x_{l,i};\mathbf{t}), \, \forall \, d=1,\ldots,k \label{eq:other-basis-functions}
\end{align}
Both Equations~\eqref{eq:0-basis-function} and \eqref{eq:other-basis-functions} are dependent on on the knot sequence $\mathbf{t} \in \mathcal{G}$ which is a sequence of non-decreasing numbers.
For a KAN, the definition of the knot vector ($\mathbf{t}$) depends on the number of grid points ($G$).
Considering the input to an activation $\phi$ is bounded between $[x_{l,i}^{\rm L}, x_{l,i}^{\rm U}]$, the knot vector $\mathbf{t}$ is defined as,
\begin{align}
    \label{eq:knot-vector}
    t_g &=
    \begin{cases}
        x_{l,i}^{\rm L} + \frac{(g - k)(x_{l,i}^{\rm U} - x_{l,i}^{\rm L})}{G}, & 0 \leq g \leq k-1, \\
        t_{g-1} + \frac{x_{l,i}^{\rm U} - x_{l,i}^{\rm L}}{G}, & k \leq g \leq G+k, \\
        x_{l,i}^{\rm U} + \frac{(g - k)(x_{l,i}^{\rm U} - x_{l,i}^{\rm L})}{G}, & G+k+1 \leq g \leq G+2k-1
    \end{cases}
\end{align}
For the input layer ($l=0$), it is straightforward to determine the bounds of the input to the network. 
For the subsequent layers, the knot vector ($\mathbf{t}$) is updated online during training \citep{Liu2024}.
Once a KAN is trained, the knot vector $\mathbf{t}$ is fixed. 
Finally, even though the defined number of grid points is $G$, the knot vector $\mathbf{t}$ comprises $G+2k$ elements with the first and last $k$ elements extending the grid points from the bounds $[x_{l,i}^{\rm L}, x_{l,i}^{\rm U}]$.
No justification is provided for the same in \citet{Liu2024}, herein we adopt their implementation.
\par
Assuming, $n_1 = \ldots = n_{L-1} = N$, a KAN is characterised by $\mathcal{O}(N^2L(G+k))$ parameters.
By contrast, for an MLP with $N$ neurons each in $L$ layers, the number of the parameters characterising the ANN is $\mathcal{O}(N^2L)$, which is lower than for KANs.
Since KANs have the ability to approximate any given input/output relationship or a function with a polynomial number of one-dimensional functions via splines, the residual rate representing the deviation between the function prediction and the data used to train the function is independent of the number of dimensions (inputs) needed to approximate the function (c.f. Theorem 2.1 in \citep{Liu2024}).
As a result, KANs typically require lower $N$ than MLPs; therefore, KANs need fewer parameters on balance to exhibit performance similar to that of an MLP.
For additional details, the reader is referred to Section 2.3 in \citet{Liu2024}.
\subsection{MIQCP formulation for B-splines}
The previous sub-section outlined the architecture of a KAN and briefly outlined how they exhibit superior efficiency in terms of parameters compared to an MLP.
Additionally, the calculations for B-splines used as activations in KANs were described in the previous section.
In this sub-section, we describe the Mixed-Integer Quadratically Constrained Programming~(MIQCP) formulation proposed by \citet{Grimstad2018} for univariate splines and adopt it for describing spline activations in KANs.
For a trained KAN,  $w^{\rm b}$, $w^{\rm s}$ and the set of coefficients $c_g$ for the spline function are fixed parameters.
Additionally, for a trained KAN the knot vectors $\mathbf{t}$ for all activations $s_{l,j,i}$ are fixed.
When embedding a KAN into an optimization formulation, the $k$-th order basis function $B_{g,k}(x_{l,i})$ is not a fixed parameter anymore and needs to be determined.
\par
The spline basis functions ($B_{g,d}(x_{l,i})$) of degree $d$ can be determined using Equations~\eqref{eq:0-basis-function} and \eqref{eq:other-basis-functions}. 
As described in Equation~\eqref{eq:0-basis-function}, zeroth-order basis functions ($B_{g,0}(x_{l,i})$) is equal to one if the input to the spline activation falls within the knot interval $[t_g, t_{g+1})$, otherwise, it is equal to zero.
For all the equations described in this sub-section we drop the indices for the layers and neurons in the spline activation function for simplicity.
\citet{Grimstad2018} introduce binary variables to represent the zeroth-order spline basis functions (Equation~\eqref{eq:zero-order-bin-var}).
\begin{equation}
    \label{eq:zero-order-bin-var}
    B_{g,0} = \{ 0,1 \}, \quad g = 0,\ldots, G+2k-1
\end{equation}
To model the discontinuity in Equation~\eqref{eq:0-basis-function} in the optimization formulation, two linear constraints (Equations~\eqref{eq:zero-order-opt-form-1} and \eqref{eq:zero-order-opt-form-2}) are introduced,
\begin{equation}
    \label{eq:zero-order-opt-form-1}
    (t_g - t_0) B_{g,0} + t_0 \leq x, \quad g = 0,\ldots, G+2k-1
\end{equation}
\begin{equation}
    \label{eq:zero-order-opt-form-2}
    (t_{g+1} - t_{G+2k}) B_{g,0} + t_{G+2k} \geq x, \quad g = 0,\ldots, G+2k-1
\end{equation}
These linear inequality constraints enforce the corresponding zeroth-order spline basis function to be active based on the value of the input belonging to the interval $t_g \leq x \leq t_{g+1}$ and are equivalent to the Big-M formulation for the discontinuous function given in Equation~\eqref{eq:0-basis-function} with the M values corresponding to the size of the intervals.
Higher-order basis functions ($B_{g,d}(x_{l,i})$) are determined by introducing Equation~\eqref{eq:other-basis-functions}, a quadratic constraint into the optimization formulation.
Following Lemma 2 in \citet{Grimstad2018}, which describes the convex combination property of B-spline basis functions, we constrain the basis functions to be non-negative (Equation~\eqref{eq:spline-basis-lb}).
\begin{equation}
    B_{g,d} \geq 0, d=0,\ldots,k
    \label{eq:spline-basis-lb}
\end{equation}
Additionally, we also consider the partition of unity cuts proposed by \citet{Grimstad2018} for all degrees $d=0,\ldots,k$ in our formulation,
\begin{equation}
    \sum_{g=0}^{G+2k-1-d} B_{g,d} = 1, \quad d=0,\ldots, k
    \label{eq:partition-of-unity-cuts}
\end{equation}
For the exact representation of a univariate B-spline, the partition of unity cut for $d=0$ is necessary.
The paritition of unity cuts for higher orders is considered to strengthen the formulation \citep{Grimstad2018}.
Finally, the k\textsuperscript{th}-order basis functions along with their respective coefficients are used to determine the output from the spline activation as described in Equation~\eqref{eq:kan-function-spline}.
Other strengthening cuts proposed in \citep{Grimstad2018} are not included in the formulation and will be discussed in Section~\ref{subsec:local_support}.
For additional details about the mathematical properties of B-spline basis functions the reader is referred to \citep{Grimstad2018}.
In the next section, we will build on this formulation for conducting the spline calculations in a KAN formulated as a MINLP.
\section{A MINLP formulation for KANs}
So far we provided a brief background on Kolmogrov Arnold Networks and introduced the MIQCP formulation proposed by \citet{Grimstad2018} for univariate B-spline functions.
In this section, we will build on the MIQCP formulation for B-splines to describe spline activation functions in KANs and propose a MINLP formulation for a fully connected KAN.
Additionally, we propose several enhancements aimed at improving the effectiveness of the MINLP formulation of the KAN.
\subsection{Description of the formulation}
\label{subsec: minlp_KAN}
Considering the optimization problems described by Equations~\eqref{eq:objective} -- \eqref{eq:bounds}, the aim is to describe a MINLP formulation for $\textsf{KAN}(\mathbf{x})$.
We consider a KAN with $L+1$ layers containing neurons defined by the vector $N := [n_0,n_1,\ldots,n_l,\ldots,n_{L-1},n_L]$, with $l=0$ and $l=L$ corresponding to the input and the output layer, respectively.
A neuron $i$ in layer $l$ is connected by an activation $\phi_{l,j,i}$ to neuron $j$ in layer $l+1$.
To model the connections in a KAN including the activations $\phi_{l,j,i}$, we introduce additional variables $\Tilde{x}_{l,j,i}$ representing the pre-activation value. 
For the input layer $l=0$, the pre-activated value $x_{0,j,i}$ is
\begin{equation}
    \Tilde{x}_{l,j,i} = x_{l,i}, \,  l=0,\ldots, L-1, \, i = 1,\ldots,n_0, \, \text{and } j=1,\ldots,n_1
    \label{eq:kan-input-enforce}
\end{equation}
The activated value for a given neuron pair ($j$,$i$) in layer $l$ is a weighted sum of the base activation function ($b$) and the spline activation function ($s$),
\begin{equation}
\begin{aligned}
    & \phi_{l,j,i}(\Tilde{x}_{l,j,i}) = w^{\rm b}_{l,j,i} b_{l,j,i}(\Tilde{x}_{l,j,i}) + w^{\rm s}_{l,j,i} s_{l,j,i}(\Tilde{x}_{l,j,i}), \\
    & l=0,\ldots,L-1, \, j=1,\ldots,n_{l+1}, \text{ and } i=1,\ldots,n_l
    \label{eq:kan-activation-enforce}
\end{aligned}
\end{equation}
Typically, \texttt{SiLU} activation is used for the base function $b$ and it is described in Equation~\eqref{eq:kan-function-base}.
For determining the spline value, the MIQCP formulation described in Equations~\eqref{eq:zero-order-bin-var} - \eqref{eq:zero-order-opt-form-2}, \eqref{eq:other-basis-functions}, \eqref{eq:partition-of-unity-cuts}, and \eqref{eq:kan-function-spline} is employed with $\Tilde{x}_{l,j,i}$ as inputs.
All the activated outputs for the layer $l$ are combined as described in Equation~\eqref{eq: output-scalar-form} and can be formulated as,
\begin{equation}
    x_{l+1,j} = \left[ \sum_{i=1}^{n_l} \phi_{l,j,i} \right] + b_{l,j}^{\rm l}, l=0,\ldots,L-1 \text{ and } j=1,\ldots,n_{l+1}
    \label{eq:kan-layer-output-enforce}
\end{equation}
While, \citet{Liu2024} do not explicitly mention the presence of a bias term ($\mathbf{b}^{\rm l}$) for a layer $l$, it is included in the implementation for sparsity regularization \citep{Liu2024b} and is a trainable parameter. 
For embedding the KAN into an optimization formulation, $\mathbf{b}^{\rm l}$ is a fixed parameter.
Thus, the description of the MINLP formulation that represents a trained KAN is complete.
It is important to note that except for Equation~\eqref{eq:kan-function-base}, all other equations in the formulation are either linear or quadratic.
In the next sub-section, we propose enhancements to strengthen the MINLP formulation of a KAN presented in this section.
\subsection{Strengthening the formulation}
In this sub-section, we propose several enhancements aimed at improving the efficiency of the MINLP formulation for a trained KAN.
\subsubsection{Feasibility-based bounds tightening (FBBT)}
\label{subsubsec: fbbt}
The training of a KAN provides us with the range of values for both inputs and outputs for a given activation $\phi_{l,j,i}$ \citep{Liu2024c}.
These range of values bound the values for $\Tilde{x}_{l,j,i}$ in the formulation.
These bounds can be used to infer bounds for the output $x_{l+1},j$ from the layer $l$ via feasibility-based bounds tightenning \citep{Puranik2017,Gleixner2017}.
Bounds can be inferred using:
\begin{equation}
\begin{aligned}
    x_{l+1,j}^{\rm L} &= \left[ \sum_{i=1}^{n_l} \phi_{l,j,i}^{\rm L} \right] + b_{l,j}^{\rm l} \\
    x_{l+1,j}^{\rm U} &= \left[ \sum_{i=1}^{n_l} \phi_{l,j,i}^{\rm U} \right] + b_{l,j}^{\rm l}
\end{aligned}
\label{eq:fbbt-layer-output}
\end{equation}
\par
Moreover, the bounds for $b$ and $s$ can be determined using the bound values for $\Tilde{x}_{l,j,i}$.
For the base activation function $b$, assuming \texttt{SiLU} activation is used, the global minimum for \texttt{SiLU} function is known.
The value of the global minimum for \texttt{SiLU} function is -0.278465 \citep{Elfwing2018} and is specified as the lower bound for $b_{l,j,i}$.
Since, \texttt{SiLU} aims to provide a smooth approximation of the \texttt{ReLU} activation function \citep{Elfwing2018}, the upper bound for $b_{l,j,i}$ is defined as (Equation~\eqref{eq:silu-ub}),
\begin{equation}
    \label{eq:silu-ub}
    b_{l,j,i}^{\rm U} = 
    \begin{cases}
        0, & \Tilde{x}_{l,j,i}^{\rm U} \leq 0 \\
        \Tilde{x}_{l,j,i}^{\rm U}, & \Tilde{x}_{l,j,i}^{\rm U} > 0
    \end{cases}
\end{equation}
\par
To bound the spline output, the bounds for the $k$-th order basis functions need to be determined. 
From Equation~\eqref{eq:spline-basis-lb}, we know that $B_{g,d}$ for any activation $l$, $j$ and $i$ is always non-negative. 
Due to the constraint $\sum_{g=0}^{G+k-1} B_{g,k} = 1$, the upper bound on $B_{g,k}$ is one since $B_{g,k}$ cannot be negative.
Hence, for a set of coefficients $\mathbf{c}$, the spline output can be bounded as,
\begin{equation}
    \begin{aligned}
        s_{l,j,i}^{\rm L} &= \sum_{g=0}^{G+k-1} \min\{ 0, c_g\} \\
        s_{l,j,i}^{\rm U} &= \sum_{g=0}^{G+k-1} \max\{ 0, c_g\}
    \end{aligned}
\end{equation}
Tighter bounds for $x_{l+1,j}$, $b_{l,j,i}$, and $s_{l,j,i}$ can be obtained using optimality-based bounds tightening techniques \citep{Puranik2017}, but is not considered in this study.
\subsubsection{Convex hull reformulation}
\label{subsubsec: hull}
In the MIQCP formulation proposed by \citet{Grimstad2018}, binary variables are introduced for each interval, and based on the value of $\Tilde{x}_{l,j,i}$ the corresponding binary variables are activated via the Big-M formulation. 
Herein, we also consider the equivalent convex hull reformulation derived from disjunctive programming for modeling Equation~\eqref{eq:0-basis-function}, which may tighten the formulation \citep{balas2018disjunctive}. Equation~\eqref{eq:0-basis-function} can be modelled as a disjunctive program,
\begin{equation}
\left[ \begin{array}{c}
B_{l,j,i,g,0}(\Tilde{x}_{l,j,i}) = 0 \\
\Tilde{x}_{l,j,i} \notin [t_g, t_{g+1}]
\end{array} \right]
\vee
\left[ \begin{array}{c}
B_{l,j,i,g,0}(\Tilde{x}_{l,j,i}) = 1 \\
\Tilde{x}_{l,j,i} \in [t_g, t_{g+1}]
\end{array} \right]
\label{eq:disjunction_0_basis}
\end{equation}
For an activation $\phi_{l,j,i}$, we introduce a set of auxiliary variables ($\mathbf{z}$) indexed over each interval in the knot vector $\mathcal{G}$, such that $z_{l,j,i,g} = B_{l,j,i,g,0} \Tilde{x}_{l,j,i}$.
To implement the convex hull reformulation, Equations~\eqref{eq:kan-ch-1} -- \eqref{eq:kan-ch-3} are introduced:
\begin{align}
    z_{l,j,i,g} & \leq t_{l,j,i,g+1} B_{l,j,i,g,0} \label{eq:kan-ch-1} \\
    z_{l,j,i,g} & \geq t_{l,j,i,g} B_{l,j,i,g,0} \label{eq:kan-ch-2} \\
    \sum_{g=0}^{G+k-1} z_{l,j,i,g} & = \Tilde{x}_{l,j,i} \label{eq:kan-ch-3}
\end{align}
Similar to the Big-M formulation, the convex hull formulation also introduces mixed-integer linear constraints to model Equation~\eqref{eq:0-basis-function}.
Convex hull reformulation is the tightest possible formulation for the disjunction in Equation~\eqref{eq:disjunction_0_basis} \citep{balas2018disjunctive}.
Whether it yields any computational benefits over the Big-M formulation will be assessed in Section \ref{sec: results} of this paper.
\subsubsection{Local support cuts}
\label{subsec:local_support}
\citet{Grimstad2018} proposed the inclusion of strengthening cuts for the MIQCP formulation for B-splines.
One class of strengthening cuts proposed by \citet{Grimstad2018} were the local support cuts which is derived from the local support property of b-splines.
The property that a basis function $B$ is bounded in $[0,1]$ is exploited by means of introducing these cuts \citep{Grimstad2018}.
Local support cuts can be appended to the formulation by the following constraint:
\begin{equation}
    B_{l,j,i,g,d} \leq \sum_{h=g}^{g+d+1-\check{d}} B_{l,j,i,h,\check{d}}, \, \check{d}=0,\ldots,k-1, \, d=\check{d}+1,\ldots,k, \, g=0,\ldots,G+k-1
    \label{eq:local_support_cuts}
\end{equation}
In Section \ref{sec: results} of this study, we will investigate the impact of appending local support cuts for optimizing over-trained KANs.
\subsubsection{Redundant cuts}
Following the idea of local support cuts, herein we propose a new class of strengthening cuts.
Again, we exploit the property that any basis function $B$ is bounded in $[0,1]$.
Considering that the constraint to determine higher order spline basis functions (Equation~\eqref{eq:other-basis-functions}) is a sum of two bilinear terms. 
The bilinear terms participating in Equation~\eqref{eq:other-basis-functions} are the previous order spline basis functions ($B_{l,j,i,g,d-1}$ and $B_{l,j,i,g+1,d-1}$) and the input to the spline $\Tilde{x}_{l,j,i}$ normalized depending on where it lies in the interval $[t_{l,j,i,g}, t_{l,j,i,g+1}]$.
Since the knot vector $\mathbf{t}$ is a sequence of non-decreasing numbers, both participating terms in the bilinear terms are also bounded in $[0,1]$, where $g=0,\ldots, G+2k-1-d$ and $d=1,\ldots,k$.
Hence, we can introduce the linear cuts:
\begin{equation}
    B_{l,j,i,g,d} \leq B_{l,j,i,g,d-1} + B_{l,j,i,g+1,d-1}, \, g = 0,\ldots, G+2k-1-d, \, d=1,\ldots,k
    \label{eq:redundant_cuts}
\end{equation}
The terms normalizing the input $\Tilde{x}_{l,j,i}$ based on where it lies in the interval $[t_{l,j,i,g}, t_{l,j,i,g+1}]$ in Equation~\eqref{eq:redundant_cuts} are assumed to be equal to one.
However, depending on the value of $\Tilde{x}_{l,j,i}$, both these terms normalizing $\Tilde{x}_{l,j,i}$ cannot take the value of one simultaneously. 
Moreover, since a basis function $B$ is bounded in $[0,1]$, Equation~\eqref{eq:redundant_cuts} is always true.
These inequalities are redundant because they do not change the solution obtained from a B-spline. 
However, they may strengthen the proposed formulation \citep{Ruiz2011}.
\par
Compared to the local support cuts defined in Equation~\eqref{eq:local_support_cuts}, fewer linear constraints are introduced to the formulation for redundant cuts.
For a KAN having two inputs trained with $G=6$ grid points using a cubic ($k=3$) B-spline, 180 redundant cuts are introduced to the formulation.
Whereas, for the same KAN, 324 local support cuts are introduced. 
Thus, local support cuts result in a formulation with higher number of constraints than redundant cuts.
Later in this paper, a more detailed analysis based on numerical experiments conducted will be presented to assess the effectiveness of these cuts.
\subsubsection{Exploiting sparsity of B-splines in KANs}
In the implementation of training a KAN, \citet{Liu2024} extend the grid by $k$ points to both left and right of the lower and upper bound in which the input to a given activation is bounded \citep{Liu2024d}.
We exploit this implementation detail to ensure input to an activation is bounded within the range determined during the training of a KAN by fixing the binary variables for extended grid points representing the zeroth-order basis functions to zero.
Consequently, the dependent higher-order basis functions determined via Equation \eqref{eq:other-basis-functions} are also fixed to zero.
\begin{equation}
    B_{\hat{d},d} = 0, \, d=0,\ldots,k \text{, } \hat{d}=0,\ldots,k-d \text{ and } \hat{d}= G+k,\ldots,G+2k-d
    \label{eq: exp_sparsity}
\end{equation}
\subsubsection{Bounding the SiLU contribution to the activation in a KAN}
\label{subsubsec: Silu}
\begin{figure}
    \centering
    \includegraphics[width=\linewidth]{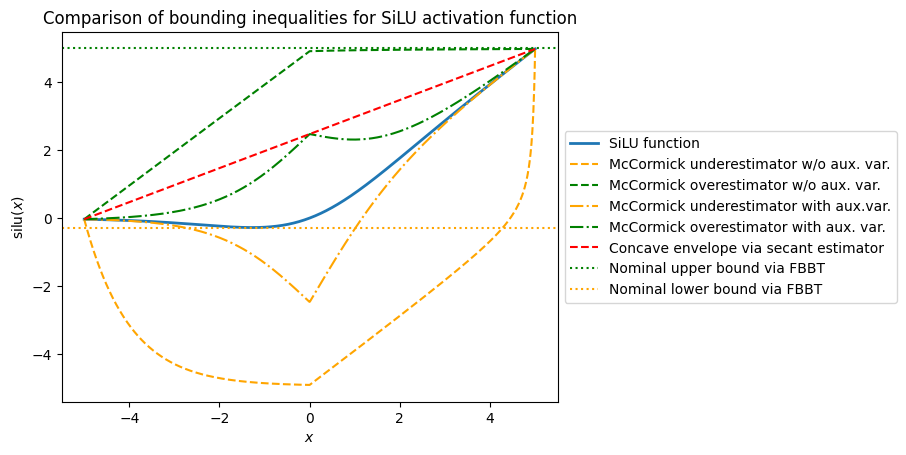}
    \caption{Comparison of different bounding strategies for the SiLU function in the range $[-5,5]$. The McCormick under and over-estimators without the introduction of an auxiliary variable are derived using \texttt{MC++} \citep{Chachuat2015}.}
    \label{fig:silu}
\end{figure}

Any activation in a KAN consists of a contribution from the trained B-spline and a base activation function (c.f. Equation \eqref{eq:kan-function}).
Typically, SiLU activation is used as the base activation function in a KAN.
Section \ref{subsubsec: fbbt} outlined how bounds were provided for SiLU contribution to an activation in KANs.
However, tighter bounds can be easily derived for the SiLU activation function.
\citet{Wilhelm2023} derived convex and concave envelopes for the SiLU activation function.
The implementation of convex and concave envelopes is more suitable within a spatial branch and bound framework where the bounds on variables are updated continuously. 
Herein, we employ the McCormick inequalities \citep{McCormick1976} to derive tighter bounds in the formulation for the SiLU activation function.
SiLU activation function can be expressed as a bilinear term where,
\begin{equation}
    b(x) = x \times y \label{eq: silu-bilinear}
\end{equation}
\begin{equation}
    y(x) = \text{sigmoid}(x) = \frac{1}{1 + \exp{(-x)}} \label{eq: sigmoid-aux-var}
\end{equation}
The range of sigmoid function is bounded in $[0, 1]$ which can be used to bound $y$. Since sigmoid function increases monotonically, tighter bounds on $y$ can be derived as $[y(x^{\rm L}), y(x^{\rm U})]$. The four linear inequalities that are added to the formulation are:
\begin{subequations}
\begin{align}
b &\geq x^{\rm L} y + x y(x^{\rm L}) - x^{\rm L} y(x^{\rm L}) \label{eq:mcc-cv1} \\
b &\geq x^{\rm U} y + x y(x^{\rm U}) - x^{\rm U} y(x^{\rm U}) \label{eq:mcc-cv2} \\
b &\leq x^{\rm U} y + x y(x^{\rm L}) - x^{\rm U} y(x^{\rm L}) \label{eq:mcc-cc1} \\
b &\leq x^{\rm L} y + x y(x^{\rm U}) - x^{\rm L} y(x^{\rm U}) \label{eq:1d}
\end{align}
\label{eq:silu_mccormick}
\end{subequations}
In Figure \ref{fig:silu} we compare the different bounding strategies for the \texttt{SiLU} activation function.
The introduction of auxiliary variable $y$ leads to tighter bounds with McCormick inequalities \citep{Najman2021} than with the nominal McCormick approach \citep{McCormick1976,Mitsos2009}, as also shown in Figure \ref{fig:silu} (c.f. the dashed green lines for bounding inequalities with an auxiliary variable for sigmoid function and the dashed yellow line for the one where no auxiliary variable is introduced).
\citet{Wilhelm2023} proposed the use of secant estimator as concave envelope to over-estimate the SiLU function.
However, the use of McCormick inequalities via the introduction of auxiliary variable $y$ for $\text{sigmoid}(x)$ leads to a tighter over-estimation of SiLU function.

\subsection{Implementation}
The MINLP formulation for a trained KAN described in Section \ref{subsec: minlp_KAN} has been implemented using \texttt{Pyomo} \citep{Bynum2021}.
A Python class named \textsf{NeuronBlockRule} is defined to create a \texttt{Pyomo} block to model a KAN activation connecting two neurons.
Equations \eqref{eq:kan-function-base}, \eqref{eq:kan-function-spline}, \eqref{eq:other-basis-functions}, and \eqref{eq:zero-order-bin-var} - \eqref{eq:kan-activation-enforce} are implemented in this class forming the \texttt{Default} formulation for a KAN.
Additionally, proposed enhancements described in Sections \ref{subsubsec: hull} -- \ref{subsubsec: Silu} are also implemented in this class. 
The proposed enhancements are activated based on the input provided by the user via an options file.
The options users must provide are described in Table \ref{tab:formulation_options}.

\begin{table}[h]
\caption{List of options that need to be passed through an options file}\label{tab:formulation_options}%
\begin{tabularx}{\textwidth}{@{}p{2cm}p{6cm}p{2cm}>{\raggedleft\arraybackslash}p{1.75cm}@{}}
\toprule
Option name & Description  & Values allowed & Default value\\
\midrule
\texttt{reformulation}      & Formulation to model the $0^{\text{th}}$-order basis functions in a B-spline in a KAN activation & \makecell[r]{ \\ Big-M, \\ Convex Hull} & Big-M\\
\texttt{exploit\_sparsity}  & Fixing basis functions corresponding to extended grid points to 0 (Equation \eqref{eq: exp_sparsity}) & \makecell[r]{0, 1} & 0  \\
\texttt{redundant\_cuts}    & Appending redundant cuts to the formulation (Equation \eqref{eq:redundant_cuts})  & \makecell[r]{0, 1} & 0  \\
\texttt{local\_support}     & Appending local support cuts to the formulation  (Equation \eqref{eq:local_support_cuts}) & \makecell[r]{0, 1} & 0 \\
\texttt{silu\_mccormick}    & Appending McCormick linear inequalities to the formulation (Equations \eqref{eq: sigmoid-aux-var}\eqref{eq:silu_mccormick}) & \makecell[r]{0, 1} & 0 \\
\botrule
\end{tabularx}
\end{table}
Hence, \texttt{Default} configuration corresponds to when all options in Table \ref{tab:formulation_options} are at their default values. \texttt{ConvexHull} configuration is when \texttt{reformulation} option is changed to ``Convex Hull" while others remain at the default values. Similarly, \texttt{ExploitSparsity} configuration is when \texttt{exploit\_sparsity} option is set to 1 and others are at their default values. \texttt{Redundant} configuration is obtained by setting the option \texttt{redundant\_cuts} to 1. On setting \texttt{local\_support} option to 1, \texttt{LocalSupport} configuration is obtained. Finally, \texttt{McCormick} configuration is obtained by setting \texttt{silu\_mccormick} option to 1.
\par
\textsf{LayerBlockRule} class inherits the \textsf{NeuronBlockRule} class to create a \texttt{Pyomo} block for a KAN layer comprising the neurons along with the inclusion of Equation \eqref{eq:kan-layer-output-enforce}.
Finally, the optimization formulation of a trained KAN is instantiated via the Python script \textsf{create\_KAN} by importing the \textsf{LayerBlockRule}. 
Additionally, the user can define the objective function and constraints relating to scaling inputs and outputs during training in this script.
\par
A Python script reads an instance of the model object resulting from training a KAN to export the values of its training parameters which include the number of layers $L+1$, the vector of number of neurons in each layer $\mathbf{N}$, coefficients to determine the spline activation value $c_{l,j,i,g}$, knot vectors $\mathbf{t}_{l,j,i}$ based on Equation~\eqref{eq:knot-vector}, the weights for the base and spline activation functions $w_{l,j,i}^{\rm b}$, and $w_{l,j,i}^{\rm s}$, respectively, biases for each layer $b_{l,j}^{\rm l}$, the lower and upper bounds on all activations $\phi_{l,j,i}$ as provided by the trained KAN model object, and the scaler parameters used to normalize the training data. 
In addition, this Python script implements one round of FBBT calculations as outlined in Section \ref{subsubsec: fbbt} to derive bounds on $x_{l+1,j}$, $b_{l,j,i}$, and $s_{l,j,i}$.
The Python script returns a JSON file with the data comprising the values of the trained parameters and bound values as described.
For instantiating a KAN, the user needs to provide the name of the JSON file in \textsf{create\_KAN} script.
\section{Results and discussion}
\label{sec: results}
Herein, we describe the computational experiments carried out to assess the effectiveness of the proposed formulation and all the enhancements.
Moreover, we discuss the results by running the six different configurations for all trained instances. 
This is followed by comparing the effect of KAN architecture on the computational effort required to globally optimize over-trained KANs. 
We discuss the impact of grid size, number of neurons, layers, and inputs. 
Finally, we compare the computational effort needed to optimize over-trained MLPs and KANs.
\subsection{Computational experiments}
\label{subsec:exp}
To understand the effect of computational effort required to optimize over a trained KAN for varying architectures, two standard test functions for optimization are considered, namely Rosenbrock and peaks function. 
Peaks function is a two-dimensional function expressed as:
\begin{align}
    f_{\rm peaks}(x_1, x_2) &= 3(1 - x_1)^2 \exp[{-x_1^2 - (x_2 + 1)^2}] - 10 \left( \frac{x_1}{5} - x_1^3 - x_2^5 \right) \exp[{-x_1^2 - x_2^2}]  \\ 
    & - \frac{\exp[{-(x_1 + 1)^2 - x_2^2}]}{3}
    \label{eq:rosenbrock}
\end{align}
The domain $D$ considered is $\{x_1, x_2 \in \mathbb{R}: -3 \leq x_1, x_2 \leq 3\}$.
In this domain the global minimum for $f_{\rm peaks}$ is -6.551 at $x_1 = 0.228$ and $x_2 = -1.626$ \citep{Schweidtmann2019}.
\par
Rosenbrock function is an $n$-dimensional function expressed as:
\begin{equation}
    f_{\rm ros}(\mathbf{x}) = \sum_{i=1}^{n-1} [100(x_{i+1} - x_i^2)^2 + (x_i - 1)^2]
\end{equation}
For Rosenbrock function, $D = \{x_1, \ldots, x_n \in  \mathbb{R}: -2.048 \leq x_1, \ldots, x_n \leq 2.048\}$, global minimum of $f_{\rm ros}$ is 0, at $x_1, \ldots, x_n = 1$ \citep{Rosenbrock1960}.
Computational experiments are conducted with the Rosenbrock function for $n_0=\{3,5,10\}$.
\par
\begin{table}[b]
\caption{Details of varying architectures used for computational experiments}\label{tab:architecture}%
\begin{tabular}{@{}llll@{}}
\toprule
Function & \# grid points  & \# neurons (\# grid points) & \# layers (\# grid points)\\
\midrule
$f_{\rm peaks}$         & $\{3,6,\ldots,24,27\}$   & $\{2,3,\ldots,9,10\} (15)$  & $\{1,2,3,4,5\} (15)$  \\
$f_{\rm ros} (n_0=3)$     & $\{3,6,\ldots,24,27\}$   & $\{2,3,\ldots,9,10\} (12)$  & $\{1,2,3,4,5\} (12)$  \\
$f_{\rm ros} (n_0=5)$     & $\{3,6,\ldots,24,27\}$   & $\{2,3,\ldots,7,8\}  (6)$   & $\{1,2,3,4,5\} (6)$  \\
$f_{\rm ros} (n_0=10)$    & $\{3,6,\ldots,24,27\}$   & $\{2,3,4,5,6\} (3)$         & $\{1,2,3,4,5\} (3)$  \\
\botrule
\end{tabular}
\end{table}
\par
Data is generated for both Rosenbrock and peaks functions, and KANs are trained. 
The details of the number of samples used for training and testing, the resulting root mean square error on both training and testing subsets, and the time needed to train the KAN are provided in the Appendix.
The effect of grid size, the number of neurons in a layer and the number of layers on the computational effort required to optimize a trained KAN is investigated.
The details of varying architectures are provided in Table \ref{tab:architecture}.
For both functions, the grid size is varied from 3 to 27, with a step increase of 3 for KAN with one layer and two neurons. 
The number of neurons in one layer is varied from 2 to 10 for a fixed grid size, and finally, the number of layers is varied from 1 to 5, with two neurons in each layer for a fixed grid size. 
The fixed grid size is reported in the corresponding column in Table \ref{tab:architecture} for investigating the effect of the number of neurons and layers.
\par
By considering different architectures, 86 instances of KANs are generated and tested on six different configurations as outlined in Table \ref{tab:formulation_options}.
All instances are solved using the global MINLP solver \texttt{SCIP} (v9.0.1) \citep{Bolusani2024} via AMPL Solver Library (ASL) in \texttt{Pyomo}. 
A default gap tolerance of zero is used for convergence of the solver \citep{AMPL2024}. 
A time limit of 2 hours is set if \texttt{SCIP} is unable to converge to the gap tolerance used. 
Additionally, multi-layer perceptrons with ReLU activation are trained using the same dataset used to train the KANs and solved using \texttt{SCIP} with the same settings by importing the trained network using \texttt{OMLT} \citep{Ceccon2022}.
A performance comparison between the computational effort needed to optimize over-trained MLPs and KANs is then carried out. 
All computational experiments are carried out on a Dell Latitude 7440 laptop with 13\textsuperscript{th} generation Intel \textsuperscript{\textregistered} Core \textsuperscript{\textregistered} i7-1365U processor @1.80 GHz with 16 GB of RAM and running Windows 10.
\subsection{Comparison of configurations}
\label{subsec:config_compare}
\begin{figure}
    \centering
    \includegraphics[width=0.65\textwidth]{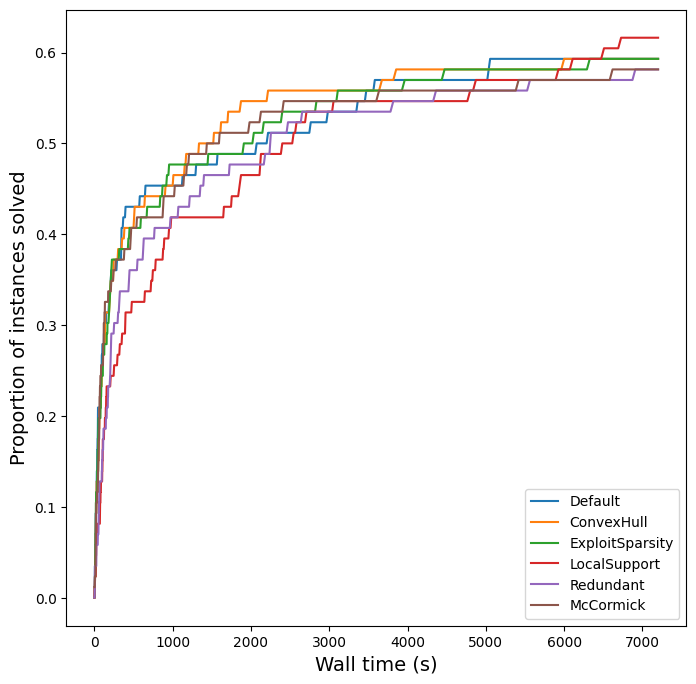}
    \caption{Performance profiles comparing the different formulation configurations of KANs}
    \label{fig:perf_profile}
\end{figure}

For an overall comparison of the different configurations considered, a modified performance profile is presented \citep{Dolan2002} in Figure \ref{fig:perf_profile}, showing the proportion of instances solved by a configuration at a given time within the time limit of 2 hours. 
A large proportion (around $33\% - 45\%$) of instances are solved within 10 minutes, depending on the configuration. 
However, regardless of the configuration, few additional instances are solved in more than 10 minutes and less than 2 hours, with around $60\%$ of the instances solved in total. 
This is a modest increase from the proportion of instances solved at 10 minutes, indicating that the computational effort required to optimize over a trained KAN increases significantly with an increase in its size (a more detailed discussion is provided in the upcoming subsections).
\par
Among different configurations, maximum instances at the end of the time limit are solved using \texttt{LocalSupport} configuration. 
However, for relatively simple instances, \texttt{LocalSupport} configuration takes longer to converge. 
For simpler instances (solved in less than 10 minutes), while local support cuts strengthen the formulation, they increase the problem's size due to additional inequalities in the formulation.
Hence, simpler instances take longer to solve with the \texttt{LocalSupport} configuration.
Similarly, due to the appending of extra constraints for the \texttt{Redundant} configuration, fewer proportion of instances are solved with this configuration at any given time. 
Moreover, the results indicate that the cuts appended in the \texttt{Redundant} configuration are relatively weaker than for \texttt{LocalSupport} configuration.
This is corroborated by the lower proportion of instances solved at the time limit using the \texttt{Redundant} configuration relative to the \texttt{LocalSupport} configuration. 
\par
Both \texttt{ExploitSparsity} and \texttt{McCormick} configurations do not offer any benefits overall, as seen in Figure \ref{fig:perf_profile}. With the \texttt{ConvexHull} configuration, a larger proportion of instances are solved within 1000 to 3000 seconds, indicating that for moderately difficult instances, the use of the convex hull over Big-M reformulation is beneficial. Overall, \texttt{Default} configuration performs the best, as corroborated by the `Shifted Geometric Wall-Time Mean' (SGWM) \citep{Mittelmann2020} reported in Table \ref{tab:SGWM}.
SGWM is calculated as,
\begin{equation}
    \text{SGWM} = \exp{\left[ \frac{\sum_{m=1}^{n_{\rm ins}} \ln(\tau_m + \alpha)}{n_{\rm ins}} \right]} - \alpha
    \label{eq: sgwm}
\end{equation}
\begin{table}[t]
\caption{SGWM for different configurations over 86 instances of KANs}\label{tab:SGWM}%
\begin{tabular}{@{}lrrr@{}}
\toprule
Configuration & \texttt{Default}  & \texttt{ConvexHull} & \texttt{ExploitSparsity} \\
\midrule
SGWM (s)    & 799.63   & 809.05  & 832.35 \\
\midrule
Configuration & \texttt{Redundant} & \texttt{LocalSupport} & \texttt{McCormick} \\
\midrule
SGWM (s)    & 1140.80 & 1235.14 & 841.75  \\
\botrule
\end{tabular}
\end{table}
In Equation \eqref{eq: sgwm}, $\tau_m$ refers to the wall time required to optimize an instance of KAN, $n_{\rm ins}$ is the number of instances and for this study 86 instances of KANs are solved, and $\alpha$ is the shift parameter, set to 5.
Henceforth, all subsequent discussions will be based on the \texttt{Default} configuration of KANs.
\subsection{Effect of grid-size}
For KANs, number of grid points used to approximate a function using B-splines is a hyperparameter. 
Bias-variance tradeoff is linked to the number of grid points used. The use of a coarse grid leads to underfitting and a fine grid leads to overfitting \citep{Liu2024}. 
In this section, we discuss the impact of varying grid size on the root mean square error (RMSE) on the testing set and the computational effort required to optimize the trained KAN. 
From Table \ref{tab:effect_grid} it is evident that on increasing the number of grid points first the test error decreases and then starts increasing.
Such behavior indicates overfitting and is consistent with the results reported in \citet{Liu2024}. 
However, this is not uniformly observed (c.f. the test RMSE for $G=12, \ldots, \, 24$) for the case of $f_{\rm peaks}$ possibly due to the stochastic nature of the training \citep{Liu2024}. 
For all the test functions, it is generally observed that on increasing the grid-size the computational effort required to optimize the KAN increases.
\begin{table}[b]
\centering
\caption{Trade-off between RMSE on the testing set and wall time (in seconds) denoted as `Time' required to solve a trained KAN with one hidden layer containing two neurons by varying the grid size $G$.}
\label{tab:effect_grid}
\begin{tabular}{r r r  r r  r r  r r}
\toprule[1.5pt]
\multirow{2}{*}{$G$} & \multicolumn{2}{c}{$f_{\rm peaks} (n_0=2)$} & \multicolumn{2}{c}{$f_{\rm ros} (n_0=3)$} & \multicolumn{2}{c}{$f_{\rm ros} (n_0=5)$} & \multicolumn{2}{c}{$f_{\rm ros} (n_0=10)$} \\ 
\cmidrule(rr){2-3} \cmidrule(rr){4-5} \cmidrule(rr){6-7} \cmidrule(rr){8-9}
 & RMSE & Time (s) & RMSE & Time (s) & RMSE & Time (s) & RMSE & Time (s) \\ 
\midrule
3  & 4.12E-01 & 5   & 1.30E-03 & 6   & 3.21E-01 & 40   & 3.71E-01 & 7200 \\
6  & 3.87E-01 & 5   & 1.19E-03 & 16  & 3.20E-01 & 46   & 3.73E-01 & 7200 \\
9  & 3.68E-01 & 12  & 1.20E-03 & 28  & 3.18E-01 & 95   & 3.74E-01 & 7200 \\
12 & 3.52E-01 & 13  & 1.20E-03 & 37  & 3.19E-01 & 7200 & 3.74E-01 & 347  \\
15 & 6.49E-01 & 25  & 1.21E-03 & 89  & 3.18E-01 & 577  & 3.73E-01 & 364  \\
18 & 5.42E-01 & 30  & 1.21E-03 & 130 & 3.17E-01 & 7200 & 3.73E-01 & 7200 \\
21 & 3.58E-01 & 27  & 1.22E-03 & 84  & 3.19E-01 & 341  & 3.70E-01 & 7200 \\
24 & 2.95E-01 & 82  & 1.23E-03 & 85  & 3.23E-01 & 5053 & 3.97E-01 & 7200 \\
27 & 4.86E-01 & 210 & 1.23E-03 & 185 & 3.49E-01 & 7200 & 4.16E-01 & 2756 \\ 
\bottomrule[1.5pt]
\end{tabular}
\end{table}
Again, some outliers are observed in Table \ref{tab:effect_grid}, see the wall time taken for the case of $f_{\rm ros}(n_0=3)$ for $G=15,\,\ldots,\,21$ for an example.
A possible reason for the occurence of outliers is that some of the trained activations may exhibit non-smooth behaviour, which adversely impacts the computational effort required to optimize a trained KAN. 
A more detailed analysis is not conducted in this work.

\subsection{Effect of number of neurons}
As in the previous sub-section a similar analysis is presented in Table \ref{tab:effect_neurons} to investigate the impact of increasing the number of neurons in a KAN.
\begin{table}[h]
\centering
\caption{Trade-off between RMSE on the testing set and wall time (in seconds) denoted as `Time' required to solve a trained KAN for all cases except $f_{\rm ros}(n_0=10)$ and the relative gap denoted as `Gap' at the end of time limit is reported for the $f_{\rm ros}(n_0=10)$ case. Trained KANs with one hidden layer containing $n_1$ neurons with grid points as defined in Table \ref{tab:architecture} are under consideration.}
\label{tab:effect_neurons}
\begin{tabular}{r r r  r r  r r  r r}
\toprule[1.5pt]
\multirow{2}{*}{$n_1$} & \multicolumn{2}{c}{$f_{\rm peaks} (n_0=2)$} & \multicolumn{2}{c}{$f_{\rm ros} (n_0=3)$} & \multicolumn{2}{c}{$f_{\rm ros} (n_0=5)$} & \multicolumn{2}{c}{$f_{\rm ros} (n_0=10)$} \\ 
\cmidrule(rr){2-3} \cmidrule(rr){4-5} \cmidrule(rr){6-7} \cmidrule(rr){8-9}
 & RMSE & Time (s) & RMSE & Time (s) & RMSE & Time (s) & RMSE & Gap (\%) \\ 
\midrule
2  & 6.49E-01 & 17    & 1.20E-03 & 36   & 3.20E-01 & 44   & 3.71E-01 & 156.90 \\
3  & 6.11E-01 & 121   & 7.14E-06 & 151  & 2.63E-01 & 7200 & 3.46E-01 & 452.02 \\
4  & 1.23E-02 & 207   & 7.79E-06 & 3459 & 6.76E-02 & 7200 & 3.19E-01 & 229.22 \\
5  & 4.12E-03 & 391   & 5.44E-05 & 5038 & 1.18E-04 & 7200 & 2.86E-01 & 354.60 \\
6  & 4.80E-03 & 1292  & 7.56E-05 & 7200 & 2.00E-04 & 7200 & 2.60E-01 & 4868.61 \\
7  & 2.81E-03 & 2066  & 9.84E-05 & 7200 & 1.08E-04 & 7200 &          &      \\
8  & 1.32E-03 & 3355  & 8.32E-05 & 7200 & 1.06E-04 & 7200 &          &      \\
9  & 1.99E-03 & 2978  & 8.38E-03 & 7200 &          &      &          &      \\
10 & 4.64E-03 & 3578  & 1.53E-04 & 7200 &          &      &          &      \\ 
\bottomrule[1.5pt]
\end{tabular}
\end{table}
Again, on increasing the number of neurons, trends of over-fitting are seen for all test functions except $f_{\rm ros} (n=5)$. 
Moreover, the computational effort required to optimize the KAN increases rapidly on increasing the number of neurons. 
In Table \ref{tab:effect_neurons} the trade-off between the wall time needed to achieve convergence or the relative gap at the time limit versus the test RMSE is markedly evident compared to the effect of grid-size. 
For each neuron added in layer $l$, $n_{l-1} + n_{l+1}$ activations are introduced in the KAN. 
For each activation, B-spline and SiLU equations need to be introduced into the formulation, thereby significantly increasing the computational effort required to solve the KAN. 
For the KAN instances considered here, the subsequent layer contains one neuron representing the output and the preceeding layer comprises all the inputs. 
Hence, for KANs with higher number of inputs, the computational effort needed to optimize increases more rapidly. 
For instance, considering $f_{\rm ros} (n=3)$ and $f_{\rm ros} (n=5)$, instances with more than 5 and 2 neurons, respectively cannot converge within the time limit. 
Remarkably, for $f_{\rm ros} (n=3)$, for the KAN with 3 neurons, one can obtain a highly accurate solution. 
Moreover, the low wall time (151 seconds) needed to solve the KAN to global optimality provides an added benefit, thereby highlighting the trememdous potential of KANs as surrogate models that can be solved with little computational effort while maintaining a high level of accuracy.
Hence, pruning of neurons \citep{Liu2024} should be considered to reduce the number of neurons while maintaining similar level of accuracy if the goal is to optimize over a trained KAN.
\subsection{Effect of layers}
Continuing with the discussion on the impact of KAN architecture, here we discuss the impact of number of layers on the trade-off between test RMSE and wall time required to optimize a trained KAN.
\begin{table}[t]
\centering
\caption{Trade-off between RMSE on the testing set and wall time (in seconds) denoted as `Time' required to solve a trained KAN for all cases except $f_{\rm ros}(n_0=10)$ and the relative gap denoted as `Gap' at the end of time limit is reported for the $f_{\rm ros}(n_0=10)$ case. Trained KANs with $L+1$ layers including the input and output layers, with each hidden layer containing two neurons with grid points as defined in Table \ref{tab:architecture}, are under consideration.}
\label{tab:effect_layers}
\begin{tabular}{r r r  r r  r r  r r}
\toprule[1.5pt]
\multirow{2}{*}{$L+1$} & \multicolumn{2}{c}{$f_{\rm peaks} (n_0=2)$} & \multicolumn{2}{c}{$f_{\rm ros} (n_0=3)$} & \multicolumn{2}{c}{$f_{\rm ros} (n_0=5)$} & \multicolumn{2}{c}{$f_{\rm ros} (n_0=10)$} \\ 
\cmidrule(rr){2-3} \cmidrule(rr){4-5} \cmidrule(rr){6-7} \cmidrule(rr){8-9}
 & RMSE & Time (s) & RMSE & Time (s) & RMSE & Time (s) & RMSE & Gap (\%) \\ 
\midrule
3  & 6.49E-01 & 17    & 9.82E-02 & 45   & 3.53E-01 & 102  & 3.71E-01 & 156.90 \\
4  & 1.80E-01 & 346   & 1.22E-03 & 1568 & 3.40E-01 & 188  & 5.55E-01 & 736.69 \\
5  & 2.65E-01 & 1119  & 6.82E-02 & 7200 & 3.22E-01 & 2208 & 3.74E-01 & 43.83 \\
6  & 1.69E-01 & 648   & 5.35E-01 & 7200 & 8.65E-01 & 7200 & 3.73E-01 & 25.26 \\
7  & 1.59E00  & 7200  & 1.33E-03 & 288  & 4.87E-01 & 7200 & 3.85E-01 & 95.12 \\
\bottomrule[1.5pt]
\end{tabular}
\end{table}
In Table \ref{tab:effect_layers}, the total number of layers $L+1$, which includes the input and output layers is reported. Hence, if the reported value for $L+1$ is three, the number of hidden layers in the KAN is one.
On varying the number of layers, there is a marked deterioration in the test RMSE for all the functions. All KANs are trained using the same set of hyperparameters and tuning them may improve the performance of KANs. It is important to note that efficient training of KANs is not the focus of this study. A discernible trade-off as seen previously is not observed in terms of computational effort on increasing number of layers. To explain the reason for this, we consider a KAN with one hidden layer and two neurons. The number of inputs is assumed to be the same as computational experiments, i.e., $2,\, 3,\, 5, \text{ and }10$. For a KAN with one hidden layer having two neurons with these number of inputs there are 6, 8, 12 and 22 activations respectively. On adding one neuron to the hidden layer, the number of activations increase to 9, 12, 18 and 33 respectively. On the other hand, if to this KAN, a hidden layer with two neurons is appended then the number of activations increase to 10, 12, 16 and 26 respectively. Thus, for KANs having fewer inputs, increasing the number of neurons in the immediate hidden layer after the input layer introduces fewer activations than adding a hidden layer with two neurons. However, for KANs with five or more inputs, adding an additional hidden layer introduces fewer activations. The computational effort required crucially depends on the number of activations in the network and the number of grid points used to define the activation. This is evident for $f_{ros}(n=5)$ and $f_{ros}(n=10)$ where a high proprtion of deeper networks can be solved for $f_{ros}(n=5)$ and the relative gaps at the end of the time limit are smaller for $f_{ros}(n=10)$ compared to a KAN with higher number of neurons. With a more careful approach to training, better test RMSEs may be achieved with deeper networks.
\subsection{Effect of number of inputs}
\begin{table}[t]
    \caption{The effect of the number of inputs on different configurations of KAN. This table reports the shifted geometric wall time mean (SGWM) in seconds and number of instances that are solved within the time limit for different cases of number of inputs. The total number of instances solved are shown in parentheses next to the number of inputs.}
    \label{tab:category}
    \centering
    \begin{tabular}{l r r r r | r r r r}
        \toprule
        Configuration & \multicolumn{4}{c|}{SGWM (s)} & \multicolumn{4}{c}{\# solved} \\
        \cmidrule(lr){2-5} \cmidrule(lr){6-9}
        & 2 & 3 & 5 & 10 & 2 (23) & 3 (23) & 5 (21) & 10 (19) \\
        \midrule
        \texttt{Default}         & 181.89 & 453.25 & 1407.98 & 4992.79 & 22 & 16 & 10 & 3 \\
        \texttt{ConvexHull}      & 207.36 & 412.83 & 1488.19 & 4743.12 & 22 & 16 & 10 & 3 \\
        \texttt{ExploitSparsity} & 213.81 & 501.14 & 1397.83 & 4414.79 & 22 & 16 & 10 & 3 \\
        \texttt{Redundant}       & 342.36 & 706.92 & 1786.42 & 5270.19 & 22 & 15 & 10 & 3 \\
        \texttt{LocalSupport}    & 399.11 & 726.47 & 1942.14 & 5540.67 & 22 & 17 & 11 & 3 \\
        \texttt{McCormick}       & 210.68 & 459.15 & 1383.96 & 5307.08 & 21 & 16 & 11 & 2 \\
        \bottomrule
    \end{tabular}
\end{table}
In the discussion of results in the previous sub-section the effect of number of inputs to the KAN is noted. 
In this section, a deeper analysis is conducted on the impact of number of inputs to the KAN on the computational effort needed to optimize them.Table \ref{tab:category} reports the proportion of instances solved and SGWM for KANs with 2, 3, 5 and 10 inputs. 
It is clearly evident that for all the configurations, the number of instances that can be solved decreases, and SGWM increases with an increase in the number of inputs. 
Similar behaviour is expected if there are multiple outputs in a KAN due to the nature of the KAN architecture. 
While the benefits of \texttt{ExploitSparsity} configuration were not clearly evident in Section \ref{subsec:config_compare}, for KANs with a large number of inputs, fixing basis functions corresponding to extended grid points resulting in smaller models leads to lower mean runtimes for KANs with five and ten inputs. 
Similarly, minor benefits are also observed for \texttt{McCormick} configuration for KANs with five inputs. The observation made earlier regarding the benefits of using convex hull formulation over Big-M formulation can also be corroborated here (c.f. SGWM for \texttt{ConvexHull} configuration for 3 and 10 inputs in Table \ref{tab:category}). 
Finally, as seen earlier, no computational benefits can be derived by using \texttt{Redundant} configuration.

\subsection{Comparison with multilayer perceptrons}
\begin{table}[b]
    \centering
    \caption{Description of architectures of different MLPs trained and how they are denoted}
    \label{tab:denote_mlp}
    \begin{tabular}{cl}
        \toprule
        MLP & Description \\
        \midrule
        1 & One hidden layer with 16 neurons \\
        2 & One hidden layer with 64 neurons \\
        3 & Two hidden layers with 16 neurons each \\
        4 & Two hidden layers with 64 neurons each \\
        5 & Three hidden layers with 16 neurons each \\
        6 & Three hidden layers with 64 neurons each \\
        \bottomrule
    \end{tabular}
\end{table}

\begin{table}[t]
\centering
\caption{Performance comparison between the best KAN and MLPs of varying architectures. The best KAN for $f_{\rm peaks}$, $f_{\rm ros}(n=3)$, $f_{\rm ros}(n=5)$, and $f_{\rm ros}(n=10)$ is one layer with seven neurons with \texttt{McCormick} configuration, one layer with four neurons with \texttt{Exploitsparsity} configuration, one layer with five neurons with \texttt{Exploitsparsity} configuration, and one layer with two neurons having 15 grid points with \texttt{ConvexHull} configuration, respectively. The description of architectures of MLPs follows Table \ref{tab:denote_mlp}. Best known indicates the value of the objective at convergence or at the time limit denoted as `PB', and best possible solution represents the dual bound at convergence or at the time limit denoted as `DB'.}
\label{tab:mlp}
\begin{tabular}{clrrrrrrr}
\toprule
\multicolumn{1}{l}{}        &               & \textsf{KAN}      & \textsf{1} & \textsf{2} & \textsf{3} & \textsf{4} & \textsf{5} & \textsf{6} \\
\midrule
\multirow{4}{*}{$f_{\rm peaks}$}    & RMSE     & 2.8E-3 & 7.3E-1 & 6.1E-1 & 5.2E-1 & 2.2E-1 & 2.7E-1 & 5.7E-2 \\
                            & Time      & 1138     & 0.01     & 2        & 2        & 604      & 5        & 7200     \\
                            & PB    & -6.551   & -2.126   & -2.844   & -3.305   & -6.293   & -5.564   & -6.775   \\
                            & DB & -6.551   & -2.126   & -2.844   & -3.305   & -6.293   & -5.564   & -21.98   \\
\midrule
\multirow{4}{*}{\makecell{$f_{\rm ros}$ \\ $(n_0=3)$}} & RMSE     & 7.8E-6 & 3.4E-1 & 2.0E-1 & 2.0E-1 & 5.3E-2 & 9.1E-2 & 4.2E-2 \\
                            & Time      & 3096     & 0.01     & 1        & 0.05     & 7200     & 2        & 7200     \\
                            & PB    & 0.00091  & -519.37  & -195.24  & -222.84  & -6.12    & -68.82   & -11.24   \\
                            & DB & 0.00091  & -519.37  & -195.24  & -222.84  & -707.5   & -68.82   & -4705.21 \\
\midrule
\multirow{4}{*}{\makecell{$f_{\rm ros}$ \\ $(n_0=5)$}} & RMSE     & 1.2E-4 & 4.9E-1 & 2.4E-1 & 2.8E-1 & 1.1E-1 & 1.7E-1 & 8.2E-2 \\
                            & Time      & 7200     & 0.01     & 2        & 1        & 7200     & 1        & 7200     \\
                            & PB    & 0.27     & -1249.75 & -724.71  & -866.48  & -276.77  & -6.38    & -102.15  \\
                            & DB & 0.19     & -1249.75 & -724.71  & -866.48  & -1834.25 & -6.38    & -9856.67 \\
\midrule
\multirow{4}{*}{\makecell{$f_{\rm ros}$ \\ $(n_0=10)$}} & RMSE & 3.7E-1 & 6.3E-1 & 3.4E-1 & 4.7E-1 & 2.6E-1 & 3.9E-1 & 2.4E-1 \\
                            & Time      & 244      & 0.01     & 3        & 1        & 7200     & 1        & 7200     \\
                            & PB    & -667.31  & -648.3   & -3156.46 & -1299.29 & -1013.9  & -544.3   & 558.89   \\
                            & DB & -667.31  & -648.3   & -3156.46 & -1299.29 & -7981.19 & -544.3   & -26727.8 \\
\bottomrule
\end{tabular}
\end{table}
For a comprehensive evaluation of the suitability of KANs as surrogate models, they are compared with MLPs. Using the same dataset, six MLPs of varying architectures are trained for each test function considered. The architectures considered are MLPs with one to three layers, each with 16 or 64 neurons. The MLPs are trained using Tensorflow \citep{tensorflow2015-whitepaper}. The trained MLPs are loaded as \texttt{Pyomo} blocks using \texttt{OMLT} \citep{Ceccon2022} and then optimized using \texttt{SCIP} \citep{Bolusani2024} on the same machine using the same settings as described in Section \ref{subsec:exp}. The results related to optimization of MLPs are reported in Table \ref{tab:mlp}.

\par
For $f_{\rm peaks}$, $f_{\rm ros}(n=3)$, and $f_{\rm ros}(n=5)$, while the smaller MLPs can be optimized very quickly, the objective obtained with those MLPs is not accurate.
The known global minimum for $f_{\rm peaks}$ is -6.551, and for  $f_{\rm ros}$ is 0. 
Relatively, on optimizing a KAN, the objective is more accurate and it can be solved in less time than larger MLPs (MLP\_2\_64 and MLP\_3\_64). 
Although, the KAN did not converge for $f_{\rm ros}(n=5)$, both the best known and best possible objectives are close to the actual value of zero, indicating the potential benefits that KANs offer as surrogate models over MLPs. 
However, for $f_{\rm ros}(n=10)$ none of the networks are able to capture the objective accurately.
In general, it is observed that MLPs do not provide accuracy with the current set of data points and hyperparameters used for training.
Using more training data or conducting hyperparameter tuning may improve the performance of MLPs, however this is not the focus of this work.
The computational effort to solve MLPs with large number of inputs does not scale as drastically as for KANs indicating that for surrogate models with more than five inputs, MLPs are a more suitable choice for surrogate models.
All test functions considered in this work have one output.
Similar behaviour is expected if more outputs were to be predicted using a KAN surrogate model since on increasing the number of outputs, the number of activations will increase similarly as that for the inputs.

\section{Conclusion}\label{sec13}
In this study, we proposed a MINLP formulation to optimize over a trained KAN. 
Additionally, we proposed some enhancements aimed at improving the effectiveness of the formulation. 
We found that the formulation based on \texttt{Default} configuration offers the best performance on balance. 
The addition of \texttt{LocalSupport} cuts, while beneficial for larger KANs, slows down the convergence for smaller KANs. 
Selective addtion of these constraints by ranking their effectiveness may further improve the effectiveness of these cuts \citep{baltean2019scoring}. 
We also observed that the formulation based on \texttt{ConvexHull} configuration also offered computational benefits for moderately difficult instances. 
The use of partition-based formulations may further improve the computational performance of the formulation \citep{kronqvist2024psplitformulationsclassintermediate}. 
Finally, providing tighter bounds for different activations via optimization-based bounds tightening \citep{Puranik2017} may further improve the effectiveness of the formulation.
\par
Additionally, we conducted a thorough investigation of the impact of the size of the KAN by varying the number of grid-points, neurons, layers and inputs on the computational effort required to optimize a KAN. We observe that the computational effort required to optimize a trained KAN crucially depends on the number of activations in the network and less strongly on the number of grid points used to define an activation. 
The number of activations introduced strongly depends on the number of inputs and outputs of the network. 
For KANs with less than five inputs, KANs offer high accuracy with relatively tractable runtimes. 
Hence, for such cases KANs offer a promising alternative as surrogate models. 
For KANs with more than five inputs, they become less attractive as surrogate models due to the prohibitive computational cost associated with optimizing them. 
For these cases, a careful consideration while deciding the architecture such as training deeper KANs rather than wider KANs can help improve tractability.

\backmatter

\section*{Declarations}
\bmhead{Funding}
All authors gratefully acknowledge financial support from Shell Global Solutions International B.V for conducting the research in this study.

\bmhead{Competing interests}
The authors have no relevant financial or non-financial interests to disclose.

\bmhead{Supplementary information}
The Zenodo repository for electronic supplementary information containing all the log files from all optimization runs can be accessed from \url{https://zenodo.org/records/14961066}.

\bmhead{Code availability}
The Pyomo formulation presented in this study can be accessed from \url{https://github.com/process-intelligence-research/optimization-over-KANs}.

\bmhead{Author contribution}
All authors contributed to the study conception and design. The formulation development, implementation, testing and analysis were performed by Tanuj Karia. The first draft of the manuscript was written by Tanuj Karia and all authors commented on previous versions of the manuscript. All authors read and approved the final manuscript. Funding was acquired by Artur M. Schweidtmann.



\noindent
If any of the sections are not relevant to your manuscript, please include the heading and write `Not applicable' for that section. 

\bigskip

\bibliography{jabref/KAN-Ref_Lib}

\end{document}